\documentclass[12pt]{amsart}

\usepackage{amsfonts,amsmath,amssymb}

\usepackage{enumerate}

\numberwithin{equation}{section}

\newtheorem{theorem}{Theorem}[section]

\newtheorem{lemma}[theorem]{Lemma}

\newtheorem{proposition}[theorem]{Proposition}

\newtheorem{example}[theorem]{Example}

\newtheorem{definition}{Definition}[section]

\newtheorem{corollary}[theorem]{Corollary}

\newtheorem{remark}[theorem]{Remark}

\newcommand{\cl}[1]{\mathcal{#1}} 

\newcommand{\bb}[1]{\mathbb{#1}}

\newcommand{\nor}[1]{\left\Vert #1\right\Vert}

\begin{document}

\title{Stable isomorphism of dual operator spaces}

\author[G.~K.~Eleftherakis]{G. K. Eleftherakis}

\address{Department of Mathematics, University of Athens,
Penepistimioupolis 157 84, Athens, Greece.}

\email{gelefth@math.uoa.gr}

\author[V.~I.~Paulsen]{V. I. Paulsen}

\address{Department of Mathematics, University of Houston,
Houston, Texas 77204-3476, U.S.A.}

\email{vern@math.uh.edu}

\author[I.~G.~Todorov]{I. G. Todorov}

\address{Department of Pure Mathematics, Queen's University Belfast, Belfast BT7 1NN, United Kingdom}

\email{i.todorov@qub.ac.uk}

\thanks{The research of the second named author was partially supported by NSF grant DMS-0600191.
The first and the third named authors were supported by EPSRC
grant D050677/1.}



\date{}

\maketitle

\begin{abstract}
We prove that two dual operator spaces $X$ and $Y$ are stably
isomorphic if and only if there exist completely isometric normal
representations $\phi$ and $\psi$ of $X$ and $Y$, respectively,
and ternary rings of operators $M_1, M_2$ such that $\phi (X)=
[M_2^*\psi (Y)M_1]^{-w^*}$ and $\psi (Y)=[M_2\phi (X)M_1^*].$ We
prove that this is equivalent to certain canonical dual operator
algebras associated with the operator spaces being stably
isomorphic. We apply these operator space results to prove that
certain dual operator algebras are stably isomorphic if and only if
they are isomorphic. We provide examples motivated by CSL algebra theory.
\end{abstract}

\section{Introduction}

K. Morita \cite{mo} developed an equivalence for rings based on their
categories of modules and proved three central theorems explaining
this equivalence relation.
A parallel Morita theory for $C^*$- and $W^*$-algebras was introduced by
Rieffel in \cite{rief}. Later Brown, Green and Rieffel \cite{bgr}
introduced the idea of stable isomorphism and proved that two
$C^*$-algebras with strictly positive elements are
strongly Morita equivalent if and only if they are stably
isomorphic in the sense that the two $C^*$-algebras obtained by
tensoring with the $C^*$-algebra of all compact operators on a
separable Hilbert space are *-isomorphic. This type of stable
isomorphism theorem is often referred to as the fourth Morita theorem,
and can often be used as an efficient way to prove some of the first
three Morita theorems. After the advent of
the theory of operator spaces and operator algebras, a parallel Morita
theory for
non-selfadjoint operator algebras was developed by Blecher, Muhly
and the second named author in \cite{bmp}. Many of the technical
results
needed to extend this theory to the setting of dual operator algebras
appear in the book of Blecher and Le Merdy \cite{bm}.
In \cite{ele1} the first named author developed a version of Morita theory for dual operator
algebras using a relation called $\Delta$-equivalence, together with a
certain category of modules over the algebras, and analogues of the first three Morita theorems were
proved. In \cite{elepaul} the first and second named authors developed the fourth part of the Morita theory, stable
isomorphism, for $\Delta$-equivalence. A different
Morita theory for
dual operator algebras has been formulated and
studied by Blecher and Kashyap \cite{bk}, \cite{kashyap}. They
have shown that it is a coarser equivalence relation than $\Delta$-equivalence,
 and have
successfully proved the first three Morita theorems in their theory.

In this paper we extend the results of \cite{ele1} and \cite{elepaul} to dual operator spaces.
We define $\Delta$-equivalence for dual operator
spaces and show that two dual operator spaces are stably
isomorphic if and only if they are $\Delta$-equivalent. Thus, we are
able to develop parts of the Morita theory in a setting where the
basic objects of study are not even rings. This
result and several of its corollaries are included in Section
\ref{s1}. We end this section by applying our results for spaces to
obtain some new results about algebras. In Section \ref{s3} we provide
examples arising from the theory of CSL algebras.

Our notation is standard. If $H$ and $K$ are Hilbert space we
denote by $H\otimes K$ their Hilbert space tensor product. For a
subset $\cl S\subseteq B(H,K)$ we denote by $\cl S'$ the commutant of $\cl S$,
by $[\cl S]$ the linear
span of $\cl S$ and by $\overline{[\cl S]}^{w^*}$ the
$w^*$-closed hull of $[\cl S]$. If $H'\subseteq H$ is a closed
subspace we let $P_{H'}$ be the orthogonal projection from $H$
onto $H'$. By $Ball(X)$ we denote the unit ball of a Banach space
$X$. For an operator algebra $A$ we denote by $pr(A)$ the set of
all projections in $A$.

Throughout the paper, we use extensively the basics of Operator
Space Theory and we refer the reader to the monographs
\cite{bm}, \cite{er}, \cite{paul} and \cite{pisier} for further details.

\section{Stably isomorphic dual operator spaces.}\label{s1}

Let $X$ be a dual operator space. A \textbf{normal representation}
of $X$ is a completely contractive $w^*$-continuous map $\phi :
X\rightarrow B(K,H)$ where $K$ and $H$ are Hibert spaces. A normal
representation $\phi : X\rightarrow B(K,H)$ is called
\textbf{non-degenerate} if $\overline{\phi (X)K}=H$ and
$\overline{\phi (X)^*H}=K$ and {\bf degenerate}, otherwise. Note
that if $\phi$ is a degenerate normal representation and if we set
$H^{\prime} = \overline{\phi(X)K},$ $K^{\prime} =
\overline{\phi(X)^*H}$ and define $\phi^{\prime}:X \to
B(K^{\prime},H^{\prime})$ by $\phi^{\prime}(x) = P_{H^{\prime}}
\phi(x) \vert_{K^{\prime}},$ then $\phi^{\prime}$ is a
non-degenerate normal representation, which we shall refer to as
the non-degenerate representation {\bf obtained from} $\phi$. If
$\phi$ is completely isometric then $\phi^{\prime}$ is completely
isometric as well. If $A$ is a unital dual operator algebra, a
\textbf{ normal representation}  of $A$ is a unital  completely
contractive $w^*$-continuous homomorphism  $\alpha : A\rightarrow
B(H)$ for some Hilbert space $H$.

If $A$ and $B$ are unital operator algebras and $X$ is an operator
space, $X$ is called an \textbf{operator $A-B$-module} if there
exist completely contractive bilinear maps $A\times X\rightarrow X$ and $X\times
B\rightarrow X$. In this case
there exist Hilbert spaces $H, K,$ completely contractive unital
homomorphisms $\pi : A\rightarrow B(H),\;\; \sigma : B\rightarrow
B(K)$ and a complete isometry $\phi : X\rightarrow B(K,H)$ such
that $\phi (axb)=\pi(a) \phi(x) \sigma (b)$ for all $a\in A, x\in
X, b\in B$ \cite[Corollary 16.10]{paul}. The triple $(\pi, \phi,
\sigma)$ is called a \textbf{CES representation} of the operator
$A-B$-module $X.$ Moreover, replacing the original $\pi$ and $\sigma$
by their direct sums with completely isometric representations, if
necessary, one may assume that $\pi$ and $\sigma$ are completely
isometric. In this case the triple $(\pi, \phi,
\sigma)$ is called a {\bf faithful CES
representation.}

If $X$ and $Y$ are dual operator spaces, we call a mapping
$\phi : X\rightarrow Y$ a {\bf dual operator space isomorphism}
if it is a surjective complete isometry which is also a $w^*$-homeomorphism.
If there exists such a mapping, we say that
$X$ and $Y$ are isomorphic dual operator spaces.
Similarly, if $A$ and $B$ are dual operator
algebras, we call a mapping $\phi : A\rightarrow B$ a {\bf dual operator algebra isomorphism}
if it is a surjective complete isometry which is also a homomorphism and a $w^*$-homeomorphism.
If there exists such a mapping, we say that
$A$ and $B$ are isomorphic dual operator algebras.

In the case that $A$ and $B$ are unital dual operator algebras and
$X$ is a dual operator space, $X$ is called a \textbf{dual
operator $A-B$-module} if it is an operator $A-B$-module and the
module actions are separately $w^*$-continuous. In this case the
triple $(\pi, \phi, \sigma )$ can be chosen with the property that
$\pi, \phi$ and $\sigma$ be $w^*$-continuous completely isometric
maps \cite[Theorem 3.8.3]{bm}. We call such a triple a {\bf
faithful normal CES representation.}

Note that since $X$ is an $A-B$-module the set $\cl C =
\left(\begin{array}{clr} A & X \\ 0 & B \end{array}\right)$ is
naturally endowed  with a product making it into an algebra and
every CES representation $(\pi, \phi, \sigma )$ as above yields a
representation $\rho: \cl C \to B(H \oplus K)$   defined by $\rho
\left( \left( \begin{array}{clr} a & x \\ 0 & b \end{array}
\right) \right) = \left( \begin{array}{clr} \pi(a) & \phi(x) \\ 0
& \sigma(b) \end{array} \right)$. When $(\pi, \phi, \sigma)$ is a
faithful CES representation, then the representation $\rho$ endows
$\cl C$ with the structure of an operator algebra. In the case $A$
and $B$ are unital $C^*$-algebras, $X$ is an operator $A-B$-module
and $(\pi, \phi, \sigma)$ is a  faithful CES representation, this
induced operator algebra structure on $\cl C$ is unique; that is,
any two faithful CES representations give rise to the same
matrix norm structures. This fact was first pointed out in
\cite[p.~11]{bp} and follows from the uniqueness of the operator
system structure on $\cl C + \cl C^*$ as can be seen from
\cite{suen} (see also \cite[3.6.1]{bm}).

In case $A$ and $B$ are $W^*$-algebras the image of the faithful
normal CES representation is $w^*$-closed and $\cl C$ can be
equipped with a dual operator algebra structure. We isolate the
following useful consequence of the above remarks.

\begin{proposition}\label{1}
Let $A_1, A_2, B_1, B_2$ be $W^*$-algebras and $X_1$ (resp. $X_2$)
be a dual $A_1-B_1$- (resp. $A_2-B_2$-) module. Let $\pi :
A_1\rightarrow A_2,\;\; \sigma : B_1\rightarrow B_2$ be normal
*-isomorphisms and $\phi : X_1\rightarrow X_2$ be a
dual operator space isomorphism which is a
bimodule map in the sense that
$$\phi (lxr)=\pi(l) \phi(x) \sigma (r), \ \ \ \  \;l\;\in A_1,\;
x\in X_1, \;r\in B_1.$$ Then the map
$$\Phi : \left(\begin{array}{clr} A_1 & X_1 \\ 0 & B_1 \end{array}\right) \rightarrow
\left(\begin{array}{clr} A_2 & X_2 \\ 0 & B_2 \end{array}\right):
 \left(\begin{array}{clr} l & x \\ 0 & r \end{array}\right)\rightarrow
 \left(\begin{array}{clr} \pi (l) & \phi (x) \\ 0 & \sigma (r) \end{array}\right)$$ is
a dual operator algebra isomorphism.
\end{proposition}

We recall some definitions from \cite{ele1} and \cite{elepaul}.
Let $I$ be a set and $\ell^2_I$ be the Hilbert space of all square
summable families indexed by $I$. Recall that if $H$ is a Hilbert
space we may identify $B(\ell^2_I \otimes H)$ with the space
$M_I(B(H))$ of all matrices of size $|I|\times |I|$ with entries
from $B(H)$ which define bounded operators on $\ell^2_I \otimes
H.$ If $X \subseteq B(H)$ is an operator space we let $M_I(X)
\subseteq M_I(B(H))$ denote the space of those operators whose
matrices have entries from $X$. We define similarly $M_{I,J}(X)$
where $I$ and $J$ are (perhaps different) index sets. In particular,
the column (resp. row) operator space $C_I(X)$ (resp. $R_I(X)$)
over $X$ is defined as $M_{I,1}(X)$ (resp. $M_{1,I}(X)$).

If $X\subseteq B(H)$ is a $w^*$-closed subspace, then it is easy
to see that $M_I(X)$ is a $w^*$-closed subspace of $M_I(B(H))$.
Moreover, if $X$ is a $w^*$-closed subalgebra of $B(H)$, then
$M_I(X)$ is a $w^*$-closed subalgebra of $M_I(B(H))$.

\begin{definition}\label{2}
(i) \cite{ele1} Let $H$ and $K$ be Hilbert spaces. Two
$w^*$-closed subalgebras $A \subseteq B(H)$ and $B \subseteq B(K)$
are called {\bf TRO-equivalent} if there exists a ternary ring of
operators (TRO) $M \subseteq B(H,K)$ such that
$A=\overline{[M^*BM]}^{w^*}$ and $B=\overline{[MAM^*]}^{w^*}.$

(ii) \cite{ele1} Two dual operator algebras $A$ and $B$ are called
{\bf $\Delta$-equivalent} if they possess completely isometric
normal representations whose images are TRO-equivalent.

(iii) \cite{elepaul} Two dual operator algebras $A$ and $B$ are
called {\bf stably isomorphic} (as algebras), if there exists a
cardinal $I$ such that the algebras $M_I(A)$ and $M_I(B)$ are
isomorphic as dual operator algebras.
\end{definition}

It is clear that stable isomorphism is an equivalence relation and
it is easy to see that the same holds for TRO-equivalence. While
it is obvious that the relation of $\Delta$-equivalence is
reflexive and symmetric, it is not apparent that it is transitive.
Nonetheless, the results of \cite{ele1} show that it is equivalent
to a certain category equivalence and hence it is also an
equivalence relation. The results of \cite{ele1} and
\cite{elepaul} show that the relations of $\Delta$-equivalence and
stable isomorphism coincide.

In this paper we generalize this result to the case of
dual operator spaces. We begin with the relevant definitions.

\begin{definition}\label{3}
(i) Let $X\subseteq B(K_1,K_2)$ and $Y\subseteq B(H_1, H_2)$ be
$w^*$-closed operator spaces. We say that $X$ is {\bf
TRO-equivalent} to $Y$ if there exist TRO's $M_1\subseteq
B(H_1,K_1)$ and $M_2\subseteq B(H_2, K_2)$ such that
$X=\overline{[M_2YM_1^*]}^{w^*}$ and
$Y=\overline{[M_2^*XM_1]}^{w^*}.$

(ii) Let $X$ and $Y$ be dual operator spaces. We say that $X$ is
{\bf $\Delta$-equivalent} to $Y$ if there exist completely
isometric normal representations $\phi$ and $\psi$ of $X$ and $Y$,
respectively, such that $\phi(X)$ is TRO-equivalent to $\psi(Y)$.

\label{4}(iii) Let $X$ and $Y$ be dual operator spaces. We say
that $X$ and $Y$ are {\bf stably isomorphic} if there exists a
cardinal $J$ and a $w^*$-continuous, completely isometric map from
$M_J(X)$ onto $M_J(Y),$ i.e., if they are isomorphic as
dual operator spaces.
\end{definition}

Blecher and Zarikian \cite[Section~6.2]{bz} define two dual operator spaces $X$ and $Y$ to be {\bf weak Morita equivalent} if $M_{I_1,J_1}(X)$ and $M_{I_2,J_2}(Y)$ are isomorphic as dual operator spaces.
Note that if $M_{I_1,J_1}(X)$ is completely isomorphic to
$M_{I_2,J_2}(Y)$ for some cardinals $I_1, I_2, J_1,J_2,$ then for a
large enough cardinal $J$ the spaces $M_J(X)$ and $M_J(Y)$ are
completely isomorphic. Thus, their definition of weak Morita equivalence is the same as our stable isomorphism. Since one goal of our research is to prove that stable isomorphism is equivalent to a type of Morita equivalence, we believe that our terminology is clearer in our context.

It is obvious that the relation of TRO-equivalence of $w^*$-closed
operator subspaces is reflexive and symmetric. We shall now prove
that it is in fact an equivalence relation. First we note that the
spaces involved can always be assumed to act non-degenerately.

\begin{proposition}\label{nondegeas}
Let $X$ and $Y$ be dual operator spaces, $\phi:X \to B(K_1,K_2),$
and $\psi: Y \to B(H_1,H_2)$ be completely isometric normal
representations with TRO-equivalent images. If $\phi^{\prime}:X
\to B(K_1^{\prime}, K_2^{\prime}),$ and $\psi^{\prime}: Y \to
B(H_1^{\prime}, H_2^{\prime})$ are the non-degenerate completely
isometric normal representations obtained from $\phi$ and $\psi,$
respectively, then the images of $\phi^{\prime}$ and
$\psi^{\prime}$ are TRO-equivalent.
\end{proposition}
\begin{proof}
Let $M_1 \subseteq B(H_1,K_1)$ and $M_2 \subseteq B(H_2,K_2)$ be
the TRO's that implement the equivalence of $\phi(X)$ and
$\psi(Y)$ and set $M_1^{\prime} = P_{K_1^{\prime}}M_1
\vert_{H_1^{\prime}}$ and $M_2^{\prime} = P_{K_2^{\prime}}M_2
\vert_{H_2^{\prime}}$. It is easy to verify that $M_1^{\prime}$
and $M_2^{\prime}$ implement a TRO-equivalence of
$\phi^{\prime}(X)$ and $\psi^{\prime}(Y).$
\end{proof}

\begin{proposition}\label{equiv}
TRO-equivalence of $w^*$-closed operator spaces is an equivalence
relation.
\end{proposition}
\begin{proof}
We need to prove that TRO-equivalence is a transitive relation.
Assume that $X\subseteq B(K_1,
K_2)$, $Y\subseteq B(H_1, H_2)$ and $Z\subseteq B(R_1, R_2)$ are
$w^*$-closed subspaces such that $X$ is TRO-equivalent to $Y$ and $Y$
is TRO-equivalent to $Z$. By Proposition \ref{nondegeas}, we may
assume that (the identity representations of) $X, Y$ and $Z$ are
non-degenerate. We fix TRO's
$$M_1\subseteq B(H_1, K_1), M_2\subseteq B(H_2, K_2), N_1\subseteq
B(H_1, R_1) \mbox{ and } N_2\subseteq B(H_2, R_2)$$ such that
$$X=\overline{[M_2YM_1^*]}^{w^*}, Y=\overline{[M_2^*XM_1]}^{w^*}, Y=\overline{[N_2^*ZN_1]}^{w^*} \mbox{ and }
Z=\overline{[N_2YN_1^*]}^{w^*}.$$ By \cite[Theorem 3.2]{ele},
there exist *-isomorphisms
$$\phi : (M_2^*M_2)^\prime\rightarrow (M_2M_2^*)^\prime \ \mbox{
and } \ \chi : (N_2^*N_2)^\prime\rightarrow (N_2N_2^*)^\prime$$
such that
$$M_2=\{T\in B(H_2,K_2): TP=\phi (P)T, \mbox{ for each } P\in pr( (M_2^*M_2)^\prime )\}$$
and
$$N_2=\{T\in B(H_2,R_2): TP=\chi (P)T, \mbox{ for each } P\in pr((N_2^*N_2)^\prime)\}.$$

Let $\cl S = pr( (M_2^*M_2)^\prime \cap  (N_2^*N_2)^\prime )$,
$$\widetilde{M_2} =\{T: TP=\phi (P)T, \mbox{ for each } P\in \cl
S\}$$
and
$$\widetilde{N_2}=\{T: TP=\chi (P)T, \mbox{ for each } P\in
\cl S\}.$$ Observe that $\widetilde{M_2} $ and $\widetilde{N_2} $
are TRO's containing $M_2$ and $N_2$, respectively. From
\cite[Lemma 2.2]{ele} it follows that

$$\overline{[\widetilde{M_2}^* \widetilde{M_2}]}^{w^*} = \cl S^\prime =
\overline{[\widetilde{N_2}^* \widetilde{N_2} ]}^{w^*}.$$ We let $
L_2 = \overline{[\widetilde{N_2}\widetilde{M_2}^*]}^{w^*}
\subseteq B(K_2,R_2).$ The space $L_2$ is a TRO since
$$\widetilde{N_2}\widetilde{M_2}^*\widetilde{M_2}
\widetilde{N_2}^*\widetilde{N_2}\widetilde{M_2}^*\subseteq
\widetilde{N_2} \cl S^\prime \cl S^\prime
\widetilde{M_2}^*\subseteq \widetilde{N_2} \widetilde{M_2}^*
\subseteq L_2.$$

Similarly, if $\cl T = pr((M_1^*M_1)^\prime \cap (N_1^*N_1)^\prime
)$ then there exist TRO's $\widetilde{M_1}\supseteq M_1,
\widetilde{N_1}\supseteq N_1 $ such that
$\overline{[\widetilde{M_1}^* \widetilde{M_1}]}^{w^*} =\cl
T^\prime = \overline{[ \widetilde{N_1}^* \widetilde{N_1}]}^{w^*}.$
As above, the space $L_1=\overline{[ \widetilde{N_1}
\widetilde{M_1}^*]}^{w^*}$ is a TRO. Since $\cl S^\prime Y \cl
T^\prime \subseteq Y$ we have
\begin{align*} \widetilde{M_2}^*\widetilde{M_2} Y
\widetilde{M_1}^*\widetilde{M_1}\subseteq Y & \Rightarrow M_2 ^*
\widetilde{M_2} Y \widetilde{M_1}^* M_1 \subseteq Y\Rightarrow
\\& M_2\stackrel{}{M_2}^*\widetilde{M_2} Y\widetilde{M_1}^*M_1M_1^*\subseteq M_2YM_1^*\subseteq X.
\end{align*}
Since $I_{K_2}\in \overline{[M_2M_2^*]}^{w^*}$ and $I_{K_1}\in
\overline{[M_1M_1^*]}^{w^*}$ we have $\widetilde{M_2}
Y\widetilde{M_1}^* \subseteq X$ and hence $X =
\overline{[\widetilde{M_2} Y\widetilde{M_1}^*]}^{w^*}$.

Similarly, we can show that
$$Y=\overline{[\widetilde{M_2}^* X\widetilde{M_1}]}^{w^*}, \ \
Z = \overline{[\widetilde{N_2} Y\widetilde{N_1}^*]}^{w^*} \ \mbox{
and } \ Y = \overline{[\widetilde{N_2}^* Z
\widetilde{N_1}]}^{w^*}.$$
Now, writing $ABC$ for $\overline{[ABC]}^{w^*}$ and $AB$ for
$\overline{[AB]}^{w^*}$ we have

$$L_2XL_1^*=\widetilde{N_2}\widetilde{M_2}^*X\widetilde{M_1}\widetilde{N_1}^*=
\widetilde{N_2}Y\widetilde{N_1}^* = Z$$ and
$$L_2^*ZL_1=\widetilde{M_2}\widetilde{N_2}^*Z\widetilde{N_1}\widetilde{M_1}^* =
\widetilde{M_2}Y\widetilde{M_1}^* = X.$$
\end{proof}

We will show later that $\Delta$-equivalence of dual operator
spaces is an equivalence relation. Note that if $A$ and $B$ are
dual operator algebras, then they could be stably isomorphic as
algebras (which requires that the map implementing the stable
isomorphism be an algebra homomorphism) or simply stably isomorphic
as dual operator spaces. However, by the operator algebra
generalization of the Banach-Stone theorem
\cite[Theorem~4.5.13]{bm} these two conditions are equivalent. In
Corollary~\ref{10} we give another proof of this fact that is
independent of the generalized Banach-Stone theorem.

We recall the following main result from \cite{elepaul}:

\begin{theorem}\label{5}
Two dual operator algebras are $\Delta$-equivalent if and only if
they are stably isomorphic as algebras.
\end{theorem}

In this section we shall generalize this result to the case of
dual operator space. Namely, we will prove the following:

\begin{theorem}\label{6}
Two dual operator spaces are $\Delta$-equivalent if and only if
they are stably isomorphic.
\end{theorem}

We now present the proof of one of the directions of Theorem
\ref{6} showing that $\Delta$-equivalence of dual operator spaces
implies stable isomorphism.

Assume, without loss of generality, that $X\subseteq B(H_1,H_2)$
and $Y\subseteq B(K_1,K_2)$ are concrete $w^*$-closed operator
spaces which are TRO-equivalent and non-degenerate. Let
$M_1\subseteq B(H_1,K_1)$ and $M_2\subseteq B(H_2,K_2)$ be
$w^*$-closed TRO's such that $\overline{[M_2 XM_1^*]}^{w^*} = Y$
and $\overline{[M_2^* YM_1]}^{w^*} = X$.

Let $$A = \left(
\begin{array}[c]{cc}
\overline{[M_2^*M_2]}^{w^*} & X\\
0 & \overline{[M_1^*M_1]}^{w^*}
\end{array}
\right) \ \ \mbox{ and } \ \ B = \left(
\begin{array}[c]{cc}
\overline{[M_2M_2^*]}^{w^*} & Y\\
0 & \overline{[M_1M_1^*]}^{w^*}
\end{array}
\right).$$ Since
$$(M_2^*M_2)X(M_1^*M_1) \subseteq M_2^*YM_1 \subseteq X,$$
the space $X$ is an
$\overline{[M_2^*M_2]}^{w^*}-\overline{[M_1^*M_1]}^{w^*}$-module
and hence $A$ is a subalgebra of $B(H_2\oplus H_1)$. Since $Y$
(resp. $X$) is non-degenerate, the relation $\overline{[M_2
XM_1^*]}^{w^*} = Y$ (resp. $\overline{[M_2^* YM_1]}^{w^*} = X$)
implies that $\overline{M_2H_2} = K_2$ (resp. $\overline{M_2^*K_2}
= H_2$). Thus, $M_2$ is non-degenerate. Taking adjoints we obtain
the relations $\overline{[M_1 X^*M_2^*]}^{w^*} = Y^*$ and
$\overline{[M_1^* Y^*M_2]}^{w^*} = X^*$ which imply that $M_1$ is
non-degenerate. It follows that the (selfadjoint) algebras
$\overline{[M_2^*M_2]}^{w^*}$ and $\overline{[M_1^*M_1]}^{w^*}$
are unital, and so $A$ is unital. One sees similarly that $B$ is a
unital $w^*$-closed subalgebra of $B(K_2\oplus K_1)$.

Let
$$M = \left(
\begin{array}[c]{cc}
M_2 & 0\\
0 & M_1
\end{array}
\right) \subseteq B(H_2\oplus H_1,K_2\oplus K_1).$$ Then $M$ is a
$w^*$-closed TRO and it is easily verified that
$$\overline{[MAM^*]}^{w^*} = B \ \ \mbox{ and } \ \ \overline{[M^*BM]}^{w^*} = A.$$ By Theorem
\ref{5}, $A$ and $B$ are stably isomorphic. Thus, there exists a
cardinal $I$ and a dual operator algebra isomorphism $\Phi
: M_I(A)\rightarrow M_I(B)$. We have that $$M_I(A) \simeq \left(
\begin{array}[c]{cc}
M_I\left(\overline{[M_2^*M_2]}^{w^*}\right) & M_I(X)\\
0 & M_I\left(\overline{[M_1^*M_1]}^{w^*}\right)
\end{array}
\right)$$ and $$M_I(B) \simeq \left(
\begin{array}[c]{cc}
M_I\left(\overline{[M_2M_2^*]}^{w^*}\right) & M_I(Y)\\
0 & M_I\left(\overline{[M_1M_1^*]}^{w^*}\right)
\end{array}
\right).$$ It is well known that $\Phi$ must carry
the diagonal of $M_I(A)$ onto the diagonal of $M_I(B)$. We claim
that $\Phi \left( \left(\smallmatrix I & 0\\ 0 &
    0\endsmallmatrix\right) \right)
= \left(\smallmatrix I & 0\\ 0 & 0\endsmallmatrix\right)$ and
$\Phi \left( \left(\smallmatrix 0 & 0\\ 0 &
I\endsmallmatrix\right) \right) = \left(\smallmatrix 0 & 0\\ 0 &
I\endsmallmatrix\right)$. To show this,  note that
$\Phi \left( \left(\smallmatrix I & 0\\ 0 &
0\endsmallmatrix\right) \right)$ is a projection in the diagonal
of $M_I(B)$ and hence there exist projections $Q$ and $P$ acting
on $K_2$ and $K_1$, respectively, such that $\Phi \left(
\left(\smallmatrix I & 0\\ 0 & 0\endsmallmatrix\right) \right) =
\left(\smallmatrix Q & 0\\ 0 & P\endsmallmatrix\right)$. Then
$\Phi \left( \left(\smallmatrix 0 & 0\\ 0 &
I\endsmallmatrix\right) \right) = \left(\smallmatrix I-Q & 0\\ 0 &
I-P\endsmallmatrix\right)$. Let $x\in M_I(X)$. Then
\begin{eqnarray*}
\Phi\left(\left(
\begin{array}[c]{cc}
0 & x\\
0 & 0
\end{array}
\right)\right) & = & \Phi\left(\left(
\begin{array}[c]{cc}
I & 0\\
0 & 0
\end{array}
\right)\left(
\begin{array}[c]{cc}
0 & x\\
0 & 0
\end{array}
\right)\left(
\begin{array}[c]{cc}
0 & 0\\
0 & I
\end{array}
\right)\right)\\
& = & \left(
\begin{array}[c]{cc}
Q & 0\\
0 & P
\end{array}
\right) \Phi\left(\left(
\begin{array}[c]{cc}
0 & x\\
0 & 0
\end{array}
\right)\right) \left(
\begin{array}[c]{cc}
I-Q & 0\\
0 & I-P
\end{array}
\right)\\
& \subseteq & \left(
\begin{array}[c]{cc}
QM_I(B(K_2))(I-Q) & QM_I(Y)(I-P)\\
0 & PM_I(B(K_1))(I-P)
\end{array}
\right).
\end{eqnarray*}
Since $\Phi$ is surjective and $Y$ is non-degenerate, it follows
that $Q = I$ and $P = 0$. The claim is proved. Since $\Phi$ is a
homomorphism, we have that $\Phi\left(\smallmatrix 0 & M_I(X)\\ 0
& 0\endsmallmatrix\right)\subseteq \left(\smallmatrix 0 & M_I(Y)\\
0 & 0\endsmallmatrix\right)$ and since $\Phi$ is onto, the last inclusion is
actually an equality. It follows that there exists a normal complete isometry
between $M_I(X)$ and $M_I(Y)$.

\bigskip

In order to prove the converse direction of Theorem \ref{6}
we need the notion of multipliers of an operator space \cite{bm}, \cite{paul}. Let $X$ be an
operator space and $M_l(X)$ be the space of all completely bounded
linear maps $u$ on $X$ for which there exist Hilbert spaces $H$ and $K$,
a complete isometry
$\iota : X\rightarrow B(H,K)$ and an operator $T\in B(K)$ such that $T \iota(X)\subseteq \iota(X)$ and
$u(x) = \iota^{-1}(T \iota(x))$, $x\in X$. Then $M_l(X)$ can be endowed
with an operator algebra structure in a canonical way and is
called the {\bf left multiplier algebra} of $X$. Similarly one
defines the right multiplier algebra $M_r(X)$ of $X$.
The operator space $X$ is an operator $M_l(X)-M_r(X)$-module; for
$l\in M_l(X)$, $r\in M_r(X)$ and $x\in X$ we write $lx = l(x)$ and $xr = r(x)$.
If $X$ is a
dual operator space then $M_l(X)$ and $M_r(X)$ are dual operator
algebras \cite[Theorem 4.7.4]{bm}. Their diagonals $A_l(X)$ and
$A_r(X)$ are thus $W^*$-algebras. Since the maps
$$A_l(X)\times X\rightarrow X: (l,x)\rightarrow lx, \;\;X\times A_r(X)\rightarrow X: (x,r)\rightarrow xr$$
are completely contractive and separately $w^*$-continuous
bilinear maps \cite[Lemma 4.7.5]{bm}, the space
\begin{equation}\label{omegaal}
\Omega (X)= \left(\begin{array}{clr} A_l(X) & X \\ 0 & A_r(X) \end{array}\right)
\end{equation}
can be canonically endowed with the structure of a dual operator
algebra (see Proposition \ref{1}).

\begin{proposition}\label{omega}
Let $X$ and $Y$ be isomorphic
dual operator spaces. Then the algebras $\Omega(X)$ and
$\Omega(Y)$ are isomorphic dual operator algebras.
\end{proposition}
\begin{proof}
Assume that $\phi : X\rightarrow Y$ is a dual operator space isomorphism. We let $\sigma :
M_l(X)\rightarrow M_l(Y)$ be given by $\sigma (u)= \phi \circ
u\circ\phi ^{-1}$. Then $\sigma$ is a completely isometric
homomorphism \cite[Proposition 4.5.12]{bm} and we can easily check
that it is $w^*$-continuous. Also, $\sigma (A_l(X))=A_l(Y)$ and
$$\phi (u x)=\phi (u(x))=\phi \circ u\circ\phi ^{-1} (\phi
(x))=\sigma (u)(\phi (x))= \sigma (u)\phi (x)$$ for all
$u\in A_l(X), x\in X$.
Similarly, the completely isometric surjection $\tau :
M_r(X)\rightarrow M_r(Y)$ given by $\tau(w)= \phi \circ w\circ\phi
^{-1}$ satisfies the identity $\phi (x w) = \phi(x)
\tau(w)$. The conclusion now follows from Proposition \ref{1}.
\end{proof}

Let $J$ be a cardinal and $X$ be a dual operator space. The
bilinear maps
$$A_l(M_J(X))\times M_J(X)\rightarrow M_J(X): (u,x)\rightarrow u(x)$$
and
$$M_J(A_l(X))\times M_J(X)\rightarrow M_J(X):$$
$$((u_{i,j}),(x_{i,j}))\rightarrow (u_{i,j})\cdot (x_{i,j}) =
\left(\sum_k u_{ik}(x_{kj})\right)_{i,j}$$ are completely
contractive and separately $w^*$-continuous by  \cite[Lemma 4.7.5,
Proposition 3.8.11]{bm}. We will need the following refinement of the statement of
\cite[Theorem~46(ii)]{bz}.

\begin{lemma}\label{cardinal}
Let $J$ be a cardinal and $X$ be an operator space. Then there
exists a *-isomorphism $\theta : M_J(A_l(X))\rightarrow
A_l(M_J(X)) $ such that
\begin{equation}\label{mod} \theta (u)(x)=u\cdot x \end{equation} for all
$u\in M_J(A_l(X)), x\in M_J(X).$ A similar statement holds for
$A_r(X)$.
\end{lemma}
\begin{proof}
We recall the operator system $S(X)=\left(\begin{array}{clr}
\mathbb{C} & X \\ X^* & \mathbb{C}
\end{array}\right)$ and its injective envelope \cite{paul}
$$I(S(X))=\left(\begin{array}{clr} I_{11}(X) & I(X) \\ I(X)^* & I_{22}(X) \end{array}\right).$$
We consider $X$ as a subspace of $I(X)$ and recall that $I(S(X))$
is a $C^*$-algebra.

For $u = (u_{ij}) \in M_J(A_l(X))$ we let $b=(b_{ij})\in
M_J(I_{11}(X))$ be the element such that $u_{ij}(x) = b_{ij}x$,
$x\in X$. Recall from \cite[Theorem~46(ii)]{bz} that there exists
a *-isomorphism $\alpha : M_J(A_l(X))\rightarrow A_l(C_J(X))$
given by $\alpha (u)(x) = b x$, for all $u\in M_J(A_l(X))$
and $x\in C_J(X)$, a *-isomorphism $\beta : A_l(C_J(X))
\rightarrow A_l(R_J(C_J(X)))$ given by
$$\beta(v)(x_i)_{i\in J}=(vx_i)_{i\in J}, \ \ \ v\in A_l(C_J(X)), (x_i)\in
R_J(C_J(X))$$ and a *-isomorphism $\gamma : A_l(R_J(C_J(X)))
\rightarrow A_l(M_J(X))$ arising from the isomorphism of $R_J(C_J(X))$ and
$M_J(X)$.

If $\theta= \gamma \circ \beta \circ \alpha :
M_J(A_l(X))\rightarrow A_l(M_J(X))$ then $\theta$ is a
*-isomorphism (and hence a $w^*$-homeomorphism) satisfying $\theta
(u)(x)=bx.$ Also,
$$\theta (u)(x)=bx= \left(\sum_{k}b_{ik}x_{kj}\right)_{i,j}
=\left(\sum_{k}u_{ik}(x_{kj})\right)_{i,j}=u\cdot x.$$
\end{proof}

In order to complete the proof of Theorem \ref{6} we will need one
more lemma.

\begin{lemma}\label{14}
Let $\cl C_X = \begin{pmatrix} B_X & X \\ 0 & A_X \end{pmatrix}$
and $\cl C_Y = \begin{pmatrix} B_Y & Y \\ 0 & A_Y \end{pmatrix}$
be concrete operator algebras acting on the Hilbert spaces $H_2
\oplus H_1$ and $K_2 \oplus K_1,$ respectively.  Suppose that
$B_X,A_X,B_Y$ and $A_Y$ are von Neumann algebras.

(i) If $\cl C_X$ and $\cl C_Y$ are TRO equivalent, then there
exist TRO's $M_1\subseteq$ $B(H_1,$ $K_1)$ and $M_2\subseteq
B(H_2,K_2)$ such that $Y = \overline{[M_2XM_1^*]}^{w^*}$ and $X =
\overline{[M_2^*YM_1]}^{w^*},$ $M_1A_XM_1^* \subseteq A_Y,
M_1^*A_YM_1 \subseteq A_X, M_2B_XM_2^* \subseteq B_Y, M_2^*B_YM_2
\subseteq B_X.$

(ii) If $\cl C_X$ and $\cl C_Y$ are $\Delta$-equivalent, then $X$
and $Y$ are $\Delta$-equivalent.
\end{lemma}
\begin{proof}
Suppose that (i) holds and assume that $\cl C_X$ and $\cl C_Y$ are
$\Delta$-equivalent.  Then there exist normal completely isometric
algebra homomorphisms, $\alpha : \cl C_X \to B(\hat{H})$ and
$\beta: \cl C_Y \to B( \hat{K})$ such that $\alpha(\cl C_X)$ and
$\beta(\cl C_Y)$ are TRO-equivalent. Note that $\alpha(\cl C_X)$
(resp. $\beta(\cl C_Y)$) has the form $\begin{pmatrix} B_Z & Z
\\ 0 & A_Z
\end{pmatrix}$ (resp. $\begin{pmatrix} B_T & T\\ 0 & A_T
\end{pmatrix}$) for a suitable decomposition $\hat{H} = \hat{H}_2\oplus\hat{H}_1$,
$\hat{K} = \hat{K}_2\oplus\hat{K}_1$, von Neumann algebras $B_Z,
A_Z, B_T, A_T$ and $w^*$-closed subspaces $Z, T$ that are
isomorphic to $X$ and $Y$, respectively, as dual operator spaces.
Thus, (ii) follows from (i).

We now prove (i). Let $P_X$ (resp. $P_Y$) denote the projection
from $H_2 \oplus H_1$ onto $H_1$ (resp. from $K_2 \oplus K_1$ onto
$K_1$). Write $\cl C_X = \cl D_X + \tilde{X},$ where $\cl D_X =
\cl C_X \cap \cl C_X^*$ and $\tilde{X} = (I-P_X)\cl C_X P_X$ is
isomorphic to $X$ as a dual operator space.
Similarly, we
decompose $\cl C_Y = \cl D_Y + \tilde{Y}.$ Let $M\subseteq
B(H_2\oplus H_1,K_2\oplus K_1)$ be a non-degenerate TRO such that
$\overline{[M\cl C_X M^*]}^{w^*} = \cl C_Y$ and $\overline{[M^*\cl
C_Y M]}^{w^*} = \cl C_X.$ Set $M_1 = P_YMP_X\subseteq B(H_1,K_1)$
and $M_2 = (I-P_Y)M(I-P_X)\subseteq B(H_2,K_2).$ Since $M$ is a
$\cl D_Y-\cl D_X$-module, we have that
$$M_1M_1^*M_1 = P_YMP_XM^*P_YMP_X \subseteq P_Y(MM^*M)P_X \subseteq M_1,$$
and hence $M_1$ is a TRO.  Similarly, we see that $M_2$ is a TRO.

Note that since $[M \cl D_XM^*]^* = [M\cl D_XM^*],$ we have that
$M\cl D_X M^* \subseteq \cl C_Y \cap \cl C_Y^* = \cl D_Y,$ and
hence $M_1A_XM_1^* \subseteq A_Y$ and $M_2B_XM_2^* \subseteq B_Y$.
Similarly, $M_1^*A_YM_1 \subseteq A_X$ and $M_2^*B_YM_2 \subseteq
B_X$.

Finally, $P_Y[M \cl D_X M^*](I-P_Y) =0,$ and it follows that
\begin{multline*}
\tilde{Y} = (I-P_Y)\cl C_Y P_Y = (I-P_Y) \overline{[M\cl D_XM^* +
M \tilde{X}M^*]}^{w^*} P_Y \\ = (I-P_Y)
\overline{[M\tilde{X}M^*]}^{w^*} P_Y = \overline{[M_2\tilde{X}
M_1^*]}^{w^*}.
\end{multline*}
Similarly, $\tilde{X} = \overline{[M_2^* \tilde{Y}M_1]} ^{w^*}$, and
hence $X$ and $Y$ are TRO-equivalent.
\end{proof}

Now we are ready to complete the proof of Theorem \ref{6}. Suppose
that $X$ and $Y$ are dual operator spaces and that there exists a
cardinal $J$ such that $M_J(X)\cong M_J(Y)$ as dual operator
spaces. We recall the unital dual operator algebras $\Omega(X)$
and $\Omega (Y)$ defined as in (\ref{omegaal}) and note that
$$ M_J(\Omega (X)) \cong \left(\begin{array}{clr} M_J(A_l(X)) & M_J(X) \\ 0 & M_J(A_r(X)) \end{array}
\right).$$ It follows from Lemma \ref{cardinal} (equation
(\ref{mod}))
and Proposition~\ref{1} that
$$\left(\begin{array}{clr} M_J(A_l(X)) & M_J(X) \\ 0 & M_J(A_r(X)) \end{array}
\right)\cong \left(\begin{array}{clr} A_l(M_J(X)) & M_J(X) \\ 0 &
A_r(M_J(X)) \end{array} \right).$$ By Proposition \ref{omega}, the
algebra on the right hand side is isomorphic as a dual operator
algebra to
$$\left(\begin{array}{clr} A_l(M_J(Y)) & M_J(Y) \\ 0 &
A_r(M_J(Y)) \end{array} \right).$$ By the same arguments, this
algebra is isomorphic to $M_J(\Omega (Y)).$ It follows from
Theorem \ref{5} that the algebras $\Omega(X), \Omega(Y) $ are
$\Delta$-equivalent as algebras.
By Lemma \ref{14} (ii), $X$ and $Y$ are $\Delta$-equivalent.
$\qquad \Box$

\bigskip

The proof of Theorem~5 is now complete. We note several immediate
corollaries.

\begin{corollary}\label{10}
If $A$ and $B$ are unital dual operator algebras then the following are equivalent:

(i) $A$ and $B$ are $\Delta$-equivalent as dual operator algebras;

(ii) $A$ and $B$ are stably isomorphic as dual operator algebras;

(iii) $A$ and $B$ are $\Delta$-equivalent as dual operator spaces;

(iv) $A$ and $B$ are stably isomorphic as dual operator spaces.
\end{corollary}

Since stable isomorphism is an equivalence relation we conclude:

\begin{corollary}
$\Delta$-equivalence of dual operator spaces is an equivalence
relation.
\end{corollary}

Let $A\subseteq B(H)$ and $B\subseteq B(K)$ be $w^*$-closed
operator algebras and $M\subseteq B(H,K)$ be a TRO such that
$$A=\overline{[M^*BM]}^{w^*} \ \mbox{ and } \ B=\overline{[MAM^*]}^{w^*}.$$
We define the $B-A$-bimodule $X \stackrel{def}{=}
\overline{[MA]}^{w^*}=\overline{[BM]}^{w^*}$ and the $A-B$-bimodule
$Y \stackrel{def}{=} \overline{[AM^*]}^{w^*}=\overline{[M^*B]}^{w^*}.$
These bimodules are important in the theory of
$\Delta$-equivalence. In \cite{ele1} they "generate" the functor
of equivalence between the categories of normal representations of
$A$ and $B.$ Also, it is proved in \cite{elepaul} that $B \simeq
X\otimes ^{\sigma h}_AY$ and $A \simeq Y\otimes^{\sigma h}_B X,$
where the tensor products are quotients of the corresponding
normal Haagerup tensor products.

\begin{corollary}\label{12}
The spaces $A,B,X,Y$ defined above are stably isomorphic.
\end{corollary}
\begin{proof} Observe that $M^*MA\mathbb{C}\subseteq A$; hence
$$M^*X\mathbb{C}\subseteq A \ \mbox{ and } \ MA\mathbb{C}\subseteq
X.$$ It follows that $X$ and $A$ are TRO-equivalent. Similarly, we
obtain that $Y$ and $B$ are TRO-equivalent. The claim now follows
from Theorem \ref{6}.
\end{proof}

In the special case of selfadjoint algebras we recapture the
following known result:

\begin{corollary}\label{13}
Let $A$ be a $W^*$-algebra and $M$ be a $w^*$-closed TRO such that
$A=\overline{[M^*M]}^{w^*}.$ Then $A$ and $M$ are stably
isomorphic.
\end{corollary}
\begin{proof} Observe that
$$MA\mathbb{C}\subseteq \overline{[MM^*M]}^{w^*}\subseteq M \ \mbox{ and } \ M^*M
\mathbb{C}\subseteq A.$$ It follows that $A$ and $M$ are
TRO-equivalent. By Theorem \ref{6}, $A$ and $M$ are stably
isomorphic.
\end{proof}

In the next result we link the $\Delta$-equivalence of two dual operator spaces $X$ and $Y$ to
that of the corresponding algebras $\Omega(X)$ and $\Omega(Y)$.

\begin{theorem}\label{9}
The dual operator spaces $X$ and $Y$ are $\Delta$-equivalent if
and only if the algebras $\Omega (X)$ and $\Omega (Y)$ are
$\Delta$-equivalent.
\end{theorem}
\begin{proof} If $X$ and $Y$ are $\Delta$-equivalent then there
exists a cardinal $I$ such that $M_I(X)$ and $M_I(Y)$ are
isomorphic as dual operator spaces. Hence, $\Omega(M_I(X))$
and $\Omega(M_I(Y))$ are isomorphic as dual operator algebras.
By Proposition \ref{1} and Lemma \ref{cardinal},
\begin{multline*}\Omega(M_I(X)) = \begin{pmatrix} A_l(M_I(X)) & M_I(X)
    \\ 0 & A_r(M_I(X)) \end{pmatrix}\\ \cong \begin{pmatrix} M_I(A_l(X))
    & M_I(X) \\0 & M_I(A_r(X)) \end{pmatrix} \cong M_I(\Omega(X)) \end{multline*}
and, similarly, $\Omega (M_I(Y)) \cong M_I(\Omega (Y))$. Thus,
$\Omega (X)$ and $\Omega (Y)$ are stably isomorphic as algebras. By
Theorem~\ref{5}, $\Omega (X)$ and $\Omega (Y)$ are
$\Delta$-equivalent.

Conversely, if $\Omega(X)$ and $\Omega(Y)$ are $\Delta$-equivalent
then, by Lemma~\ref{14} (ii), $X$ and $Y$ are $\Delta$-equivalent.
\end{proof}

\begin{theorem}\label{newbim} Let $X$ and $Y$ be $\Delta$-equivalent dual operator spaces.
If $(\pi, \phi, \sigma)$ is a normal CES representation
of the dual operator $A_l(X)-A_r(X)$-module $X$ and $\phi$ is a complete isometry,
then there exists a normal completely isometric
representation $\psi $ of $Y$ such that $\phi (X)$ is
TRO-equivalent to $\psi (Y).$
\end{theorem}
\begin{proof}
The CES triple $(\pi, \phi, \sigma)$ defines a normal
representation $\Phi$ of the algebra $\Omega
(X).$ If $l\in A_l(X)$ with $\pi (l)=0$ then $\phi
(lx)=0$ and hence $lx=0$ for all $x\in X.$ This implies that $l=0,$ and so $\pi$ is one-to-one.
Similarly $\sigma$ is one-to-one. Thus, $(\pi, \phi, \sigma)$ is a
faithful CES representation and induces the unique operator algebra
structure on $\Omega(X).$ Thus, $\Phi$ is a normal completely isometric
representation of the dual operator algebra $\Omega(X).$

By Theorem \ref{9}, $\Omega (X)$ and $\Omega
(Y)$ are $\Delta$-equivalent; by \cite[Theorem 2.7]{ele2}, there
exists a normal completely isometric representation $\Psi $ of
$\Omega (Y)$ such that $\Phi (\Omega (X))$ is TRO-equivalent to
$\Psi (\Omega (Y)).$

Let $\psi $ be the restriction of $\Psi $ to $Y.$ By Lemma
\ref{14} (i), the spaces $\phi (X)$ and $\psi (Y)$ are TRO-equivalent.
\end{proof}

\bigskip

By \cite{ele1}, $\Delta$-equivalence for dual operator algebras can be
equivalently defined in terms of a special type of isomorphism
between certain categories of representations of the algebras.
These types of category isomorphisms are in the spirit of Morita
equivalence. Thus, one would like to claim that the
representations of $\Omega(X)$ and of $\Omega(Y)$ define certain
special families of representations of $X$ and $Y$ such that $X$
and $Y$ are stably isomorphic if and only if these classes of
representations are isomorphic.  Unfortunately, the correspondence
between representations of $\Omega(X)$ and representations of $X$
is not one-to-one.

We finish this section with some applications of the above theorems.

\begin{definition} An operator space $X$ is called {\bf rigid} if
  $M_l(X) = M_r(X) = \bb C$ and {\bf *-rigid} if $A_l(X) = A_r(X) = \bb C.$
\end{definition}

Note that if $X$ is rigid, then it is *-rigid. There are many
examples of rigid and *-rigid operator spaces. For example, the
spaces $MAX(\ell^n_1)$ by a result of Zhang \cite{zhang} (see also
\cite[Exercise~14.3]{paul}) can be identified with the subspace of
the full group $C^*$-algebra of the free group on $n-1$
generators, $C^*(\bb F_{n-1}),$ spanned by the identity and the
$n-1$ generators. Moreover, Zhang argues that $I(MAX(\ell^n_1)) =
I(C^*(\bb F_{n-1}))$ and since $C^*(\bb F_{n-1})$ is a
$C^*$-subalgebra of its injective envelope it  follows from
\cite{bp} that any left multiplier of $MAX(\ell^n_1)$ necessarily
belongs to this subspace and multiplies the subspace back into
itself in the usual product. This forces the multiplier to be a
multiple of the identity. A similar argument applies for right
multipliers. Thus, $MAX(\ell^n_1)$ is rigid.

This argument given in the previous paragraph applies equally well to any subspace $X$ of a
unital $C^*$-algebra $A$ which contains the identity and for which
$I(X)=I(A).$ In this case, the left (and right) multipliers are
simply the elements of the subspace that leave the subspace
invariant under the algebra multiplication, and so it is often
quite easy to determine whether $X$ is rigid or *-rigid.

\begin{theorem}\label{riggid}
Let $X$ and $Y$ be *-rigid dual operator spaces. Then $X$ and $Y$
are stably isomorphic if and only if they are isomorphic as dual operator spaces.
\end{theorem}
\begin{proof}
If $X$ and $Y$ are stably isomorphic, then they are
$\Delta$-equivalent. Hence, by Theorem \ref{9}, $\Omega(X)$ and
$\Omega(Y)$ have completely isometric representations whose images
are TRO-equivalent. The images of these representations are two
concrete operator algebras $\cl C_X$ and $\cl C_Y$ of the type
considered in Lemma \ref{14}, with $A_X, B_X, A_Y$ and $B_Y$ all
scalar multiplies of the identity and $X$ and $Y$ replaced by
images of normal completely isometric representations, say
$\phi(X)$ and $\psi(Y).$ Hence, by Lemma \ref{14} (i), the TRO's
$M_1$ and $M_2$ arising in its proof satisfy $M_1^*M_1= M_1M_1^* =
M_2M_2^* = M_2M_2^* = \bb C.$

Now it readily follows that the spaces $M_1$ and $M_2$ are each
the span of a single unitary.  Let $M_i = \bb C U_i, i=1,2$ for
some unitaries $U_1$ and $U_2$.  Applying Lemma \ref{14}
again, we see that $\psi(Y) = U_2^* \phi(X)U_1$ and the claim
follows.
\end{proof}

\begin{corollary}  Let $A$ and $B$ be dual operator algebras for which
  $A \cap A^* = B \cap B^* = \bb C.$ Then $A$ and $B$ are stably
  isomorphic as operator spaces if and only if they are isomorphic as dual
  operator algebras.
\end{corollary}
\begin{proof} Since $B$ is a unital algebra, we have that
  $M_l(B) = M_r(B) = B$ and hence, $A_l(B) = A_r(B) = B \cap B^* = \bb
  C.$  Hence, $B,$ and similarly $A,$ is a $*$-rigid operator space.
Thus, by Theorem \ref{riggid}, $A$ and $B$ are stably
isomorphic if and only if they are isomorphic as dual operator spaces.
By the generalized Banach-Stone
theorem \cite[Thereom~3.8.3]{bm}, $A$ and $B$ are isomorphic as dual operator algebras.
\end{proof}

It is interesting to note that the hypotheses and conclusions of the above corollary
are really special to non-selfadjoint operator algebras.
In fact, we now turn our attention to a special family of
non-selfadjoint operator algebras to which our theory applies.

\begin{definition}
Let $G \subseteq \bb C^n$ be a bounded, connected, open set, i.e., a
complex domain, and let
$H^{\infty}(G) \subseteq L^{\infty}(G)$ denote the dual operator
algebra of bounded analytic functions on $G.$  We shall call $G$ {\bf
  holomorphically complete} if every weak*-continuous multiplicative
linear functional on $H^{\infty}(G)$ is given by evaluation at
some point in $G.$
\end{definition}

Recall that two complex domains $G_i \subseteq \bb C^{n_i},
i=1,2$ are called {\bf biholomorphically equivalent} if there exists a
holomorphic homeomorphism, $\varphi: G_1 \to G_2$ whose inverse is also holomorphic.

\begin{corollary} Let $G_i, i=1,2$ be complex domains that are
  holomorphically complete.  Then the following are equivalent:

(i) $G_1$ and $G_2$ are biholomorphically equivalent,

(ii) $H^{\infty}(G_1)$ and $H^{\infty}(G_2)$ are isometrically
weak*-isomorphic algebras,

(iii) $H^{\infty}(G_1)$ and $H^{\infty}(G_2)$ are isometrically
weak*-isomorphic dual Banach spaces,

(iv) $H^{\infty}(G_1)$ and $H^{\infty}(G_2)$ are stably isomorphic
dual operator spaces,

(v) $H^{\infty}(G_1)$ and $H^{\infty}(G_2)$ are $\Delta$-equivalent
dual operator algebras,

(vi) $H^{\infty}(G_1)$ and $H^{\infty}(G_2)$ are $\Delta$-equivalent
dual operator spaces.

\end{corollary}
\begin{proof}  Since $G_1$ and $G_2$ are connected sets, we have that
$H^{\infty}(G_i) \cap H^{\infty}(G_i)^* = \bb C, i=1,2.$ Also, since
these algebras are subalgebras of commutative
  C*-algebras, every contractive map between them is automatically
  completely contractive. Thus, the equivalence of (ii)--(vi) follows
  from the previous results.

Given a biholomorphic map $\varphi: G_1 \to G_2,$ composition with $\phi$
defines the weak*-continuous isometric isomorphism between the
algebras. Thus, (i) implies (ii).

Conversely, given a weak*-continuous isometric algebra isomorphism,
$\pi: H^{\infty}(G_1) \to H^{\infty}(G_2),$ let $w \in G_2,$ and let
$E_w: H^{\infty}(G_2) \to \bb C$ denote the weak*-continuous, multiplicative linear
functional given by evaluation at $w.$ Then $E_w \circ \pi:
H^{\infty}(G_1) \to \bb C$ is a weak*-continuous, multiplicative
linear functional and hence is equal to $E_z$ for some $z \in G_1.$ If we assume that $G_1 \subseteq \bb C^n,$ let $z_1, \ldots, z_n$ denote the coordinate
functions on $G_1$ and set $\varphi_i = \pi(z_i) \in H^{\infty}(G_2),$
then it readily follows that $\varphi = (\varphi_1, \ldots,
\varphi_n): G_2 \to \bb C^n$ satisfies, $\varphi(w) = z.$  Hence,
$\varphi: G_2 \to G_1.$ A similar argument with the inverse of $\pi$ shows that
$\varphi$ is a biholomorphic equivalence.  Thus, (ii) implies (i).

\end{proof}

Recalling that $\Delta$-equivalence is originally defined in terms of
a Morita-type equivalence of categories, we see that the equivalence
of (i) and (v) shows that two domains have ``equivalent'' categories
of representations in this sense if and only if they are
biholomorphically equivalent.

\section{Applications and examples}\label{s3}

In this section we prove that whenever two dual operator algebras
$A$ and $B$ are $\Delta$-equivalent, there exists a dual operator
space $X$ such that $A$ is completely isometrically isomorphic to
$M_l(X)$ and $B$ is completely isometrically isomorphic to
$M_r(X).$ We then give an example of a dual operator space $Y$ for
which $M_l(Y)$ and $M_r(Y)$ are not stably isomorphic and hence
not $\Delta$-equivalent. We also give some examples which emphasize
the difference between dual operator spaces arising from
non-synthetic CSL algebras and those arising from synthetic ones.

Let $A\subseteq B(H)$ be a unital $w^*$-closed algebra, $\Delta(A)
= A \cap A^*$ be its diagonal and $M\subseteq B(K,H)$ be a
non-degenerate TRO such that $MM^*\subseteq A.$ We call the space
$X = \overline{[AM]}^{w^*}$ the \textbf{$M$-generated $A$-module.}
In this section we fix $A$ and $M$ as above and we investigate
some properties of $X.$ Since $MM^*\subseteq A$ the space $B =
\overline{[M^*AM]}^{w^*}\subseteq B(K)$ is a unital algebra and
$XB\subseteq X$. Note that if we set $Y= \overline{[M^*A]}^{w^*},$
then $A,X,Y,B$ form the four corners of what could potentially be
a "linking" algebra of a Morita context.
For this reason, we shall call $(A,M,X,B)$ a {\bf generating tuple.}

\begin{theorem}\label{bim1} Let $(A,M,X,B)$ be a generating tuple.

(i) $M_l(X)$ is isomorphic as a dual operator algebra to $A$ and
$M_r(X)$ is isomorphic as a dual operator algebra to $B.$

(ii) The algebra $\Omega (X)$ is isomorphic as a dual operator
algebra to the algebra $D(X)=\left(\begin{array}{clr} \Delta (A) &
X
\\ 0 & \Delta (B)
 \end{array}\right).$

\end{theorem}

\begin{proof}
If $a\in Ball(A)$ we define a map $\lambda (a): X\rightarrow X$ by
letting $\lambda (a)(x)=ax.$ If $x,y\in X$ then
$$\nor{(ax, y)^t}=\nor{ \left(\begin{array}{clr}a & 0 \\ 0 & I\end{array}\right) (x, y)^t}\leq
\|a\|\|(x, y)^t\|,$$ and hence the map
$$\tau :
C_2(X)\rightarrow C_2(X), \ \ \tau ((x,y)^t)=(ax,y)^t$$ is
contractive. Similarly one can show that $\tau$ is completely
contractive. By \cite[Theorem 4.5.2]{bm}, $\lambda (a)\in
Ball(M_l(X))$.

It follows that the map $$\lambda : A\rightarrow M_l(X):
a\rightarrow \lambda (a)$$ is contractive. It is also
$w^*$-continuous by \cite[Theorem 4.7.4]{bm}. We now prove that
$\lambda $ is an isometric surjection. Using analogous arguments,
we can show that the map
$$\rho : B\rightarrow M_r(X), \ \ \  \rho (b)(x)=xb$$
is $w^*$-continuous and contractive. Let $u$ be in $M_l(X).$ By
\cite[Lemma 8.5.23]{bm} there exists a family $(m_i)_{i\in
I}\subseteq M$ of partial isometries such that $m_im_i^* \bot
m_jm_j^*$ for $i\neq j$ and $I_H=\sum_{i\in I}m_im_i^*$, the series converging
in the strong operator topology. Let $x\in
X, \xi \in K$ and $F\subseteq I$ be finite. Since the operators on
$X$ from $M_l(X)$ commute with those from $M_r(X)$, we have
\begin{align*} \sum_{i\in F}u(m_i)m_i^*x(\xi ) =& \sum_{i\in F}\rho (m_i^*x)(u(m_i))(\xi
)=\sum_{i\in F}u(\rho (m_i^*x)m_i)(\xi )\\
=& \sum_{i\in F}u(m_im_i^*x)(\xi )=u\left(\sum_{i \in
F}m_im_i^*x\right)(\xi ).\end{align*}

Since $u$ is $w^*$-continuous \cite[Theorem 4.7.1]{bm} we have
that
\begin{equation}\label{limi}
\lim_F \sum_{i\in F}u(m_i)m_i^*x(\xi ) =u(x)(\xi), \ \ \  \xi \in
K.
\end{equation}
Observe that if $F=\{i_1,...,i_n\}\subseteq I$ then
\begin{align*} \nor{ \sum_{i\in F}u(m_i)m_i^*}=&\nor{u( (m_{i_1},...,m_{i_n}) )
(m_{i_1^*},...,m_{i_n}^*)^t }\leq \\&
\|u\|_{M_l(X)}\nor{(m_{i_1},...,m_{i_n}) }
\nor{(m_{i_1^*},...,m_{i_n}^*)^t }\leq \|u\|_{M_l(X)}.\end{align*}
Hence, the net $( \sum_{i\in F}u(m_i)m_i^* )_F$ is bounded. Since
$X$ is non-degenerate the limit of the net $(\sum_{i\in
F}u(m_i)m_i^*(\xi))_F $ exists for all $\xi \in H.$ We let $a =
\sum_{i\in I}u(m_i)m_i^* \in A$, the series converging in the
strong operator topology. Observe that
\begin{equation}\label{exx}
\|a\|\leq \|u\|_{M_l(X)}.
\end{equation}
By (\ref{limi}), $ax=u(x)$ for all $x\in X$ and so $u=\lambda
(a).$ We proved that $\lambda $ is onto. By standard arguments,
equation (\ref{exx}) implies that $\lambda $ is isometric.

Let
$n\in \mathbb{N}$ and $N = \underbrace{M\oplus \dots\oplus M}_n$.
Then the $N$-generated $M_n(A)$-module is equal to
$M_n(X)=\overline{[M_n(A)N]}^{w^*}.$ By the arguments above, the
map
$$\sigma : M_n(A)\rightarrow M_l(M_n(X)): \sigma (a)(x)=ax$$ is a surjective isometry.

We recall the surjective isometry \cite[Theorem 5.10.1]{bz}
$$L: M_n(M_l(X))\rightarrow M_l(M_n(X)):
L((u_{ij})_{i,j})((x_{ij})_{i,j}) =
\left(\sum_{k}u_{ik}(x_{kj})\right)_{i,j}.$$ Since $\lambda^{(n)}
= L^{-1}\circ \sigma : M_n(A)\rightarrow M_n(M_l(X)),$ we have
that $\lambda $ is $n$-isometric. We have thus shown that $\lambda$ is a completely isometry.
Similarly, we can prove that $\rho $ is
completely isometric and surjective. By Proposition \ref{1}, the map
$$\Phi : D(X)\rightarrow \Omega (X): \left(\begin{array}{clr}a & x \\ 0 & b\end{array}\right)
\rightarrow \left(\begin{array}{clr}\lambda (a) & x \\ 0 & \rho
(b)\end{array}\right)$$ is a dual operator algebra isomorphism.
\end{proof}

\begin{corollary}\label{bim2}
If $C$ and $D$ are $\Delta$-equivalent unital dual operator
algebras then there exists a dual operator space $X$ such that
$C\cong M_l(X)$ and $D\cong M_r(X)$ as dual operator algebras.
\end{corollary}
\begin{proof}
The algebras $C$ and $D$ have completely isometric normal
representations which are TRO-equivalent. Letting $A$ be the image of $C,$ letting $M$ be the TRO that induces the equivalence and applying
Theorem \ref{bim1} to the corresponding generating tuple completes the proof.
\end{proof}

\begin{remark}\label{bim10}
{\rm The converse of Corollary \ref{bim2} does not hold. Example
\ref{nestex} shows that there exists a dual operator space $Y$
such that $M_l(Y)$ and $M_r(Y)$ are not stably isomorphic.}
\end{remark}

\begin{proposition}\label{bim3} Let $(A,M,X,B)$ be a generating tuple.
If $Y$ is a dual operator space which is $\Delta$-equivalent to the
dual operator space $X,$ then there exists a normal completely isometric
representation $\psi $ of $Y$ such that $X$ is
TRO-equivalent to $\psi (Y).$
\end{proposition}
\begin{proof}
By Theorem \ref{bim1}, $\Delta (A)$ is isomorphic to $A_l(X),$
and
$\Delta (B)$ is isomorphic to $A_r(X).$ Thus, there is a normal CES
representation of the form $(\pi, id_X, \sigma)$ of the dual
operator $A_l(X)-A_r(X)$-module $X$. Now apply Theorem \ref{newbim}.
\end{proof}

\bigskip

We recall some definitions and concepts that we will need in the
rest of the paper, see \cite{dav}. A \textbf{commutative subspace
lattice (CSL)} is a strongly closed projection lattice $\cl{L}$
whose elements mutually commute. A \textbf{CSL algebra} is the
algebra $\mathrm{Alg}\cl{L}$ of operators leaving invariant all
projections belonging to a CSL $\cl L$. In the special case where
$\cl L$ is totally ordered we call $\cl L$ a \textbf{nest} and the
algebra $\mathrm{Alg}\cl{L}$ a \textbf{nest algebra}.
There exists a smallest $w^*$-closed algebra contained in $\cl{A}$
which contains the diagonal $\Delta (\cl{A})$ of $\cl A$ and whose
reflexive hull is $\cl{A}$ \cite{a}, \cite{os1} (for the definition of
the reflexive hull of an operator algebra see \cite{dav}).
We denote this
algebra by $\cl{A}_{min}.$  If $\cl{A}= \cl{A}_{min}$ the
CSL algebra is called $\cl{A}$ \textbf{synthetic}.

\begin{proposition}\label{bim4} Let $A$ and $D$ be CSL
algebras and $M$ and $N$ be TRO's such that $MM^* \subseteq A$ and
$NN^*\subseteq D$. Set $X = \overline{[AM]}^{w^*}$ and
$Y=\overline{[DN]}^{w^*}$. Then $X$ and $Y$ are
$\Delta$-equivalent if and only if they are TRO-equivalent.
\end{proposition}
\begin{proof}
Suppose that $X$ and $Y$ are $\Delta$-equivalent. Since
$\Omega(X)$ and $\Omega(Y)$ are $\Delta$-equivalent, the algebras
$D(X)$ and $D(Y)$ defined as in Theorem \ref{bim1} (ii) are
$\Delta$-equivalent. It is easy to see that $D(X)$ and $D(Y)$ are
CSL algebras. It follows from \cite[Theorem 3.2]{ele2} that $D(X)$
and $D(Y)$ are TRO-equivalent. By Lemma \ref{14} (i), $X$ and $Y$ are
TRO-equivalent.
\end{proof}

\begin{example}\label{bim5}
\em{We now give an example of spaces which are not
$\Delta$-equivalent. Let $A$ be a CSL algebra, $B$ be a
non-synthetic, separably acting CSL algebra and $M$ and $N$ be
TRO's such that $MM^*\subseteq A$ and $NN^*\subseteq B$. Then the
spaces $X=\overline{[AM]}^{w^*}$ and
$Y=\overline{[B_{min}N]}^{w^*}$ are not $\Delta$-equivalent.
Indeed, if they were, they would be stably isomorphic. On the
other hand, Corollary \ref{12} implies that $X$ is stably
isomorphic to $A$ and $Y$ is stably isomorphic to $B_{min}.$ Thus,
the algebras $A$ and $B_{min}$ would be stably isomorphic, hence
$\Delta$-equivalent. This contradicts \cite[Theorem
3.4]{ele2}. }
\end{example}

\bigskip

Let $\cl{N}_1$ and $\cl{N}_2$ be nests acting on separable Hilbert
spaces $H_1$ and $H_2$, respectively. Recall \cite{dav} that $\cl
N_1$ and $\cl N_2$ are called similar if there exists an
invertible operator $y : H_1\rightarrow H_2$ such that $\cl N_2 =
\{yn(H_1) : n\in \cl N_1\}$. In this case there exists an order
isomorphism $\theta : \cl{N}_1\rightarrow \cl{N}_2$ which
preserves the dimension of the atoms of $\cl N_1$ and $\cl N_2$, namely,
$\theta(n)$ can be taken to be equal to the projection onto $yn(H_1)$, for all
$n\in \cl{N}_1$.
We say that the invertible operator $y\in B(H_1, H_2)$ implements
$\theta$. Let
$$Y=\{y\in B(H_1, H_2): (I-\theta(n))yn = 0, \;\;\forall \;\;n\in \cl{N}_1\}$$
and
$$Z=\{x\in B(H_2, H_1): (I-n)x\theta(n) = 0, \;\;\forall \;\;n\in \cl{N}_1\}.$$
If $C=\mathrm{Alg}\cl{N}_1$ and $D=\mathrm{Alg}\cl{N}_2$ one can
easily verify that
$$C=\overline{[ZY]}^{w^*},\;\; D=\overline{[YZ]}^{w^*},\;\; CZD\subseteq Z \ \mbox{ and } \ DYC\subseteq Y.$$

We will need the Similarity Theorem \cite[Theorem 13.20]{dav}:

\begin{theorem}\label{bim11}
For every $\delta >0$ there exists an invertible operator $y$
which implements $\theta$ such that $\|y\|<1+\delta$ and
$\nor{y^{-1}}<1+\delta$.
\end{theorem}

\begin{theorem}\label{bim12}
(i) $M_l(Z)\cong C,\;\; M_r(Z)\cong D$ as dual operator algebras.

(ii) The algebra $\Omega (Z)$ is isomorphic as a dual operator
algebra to the algebra $\left(\begin{array}{clr} \Delta (C) & Z
\\ 0 & \Delta (D)
\end{array}\right).$
\end{theorem}

\begin{proof}
We can easily check that the map
$$\tau : C_2(Z)\rightarrow C_2(Z): (x_1, x_2)^t\rightarrow (ax_1, x_2)^t$$ is completely contractive
for all $a\in C.$ So by \cite[Theorem 4.5.2]{bm} the linear map
$\lambda : C\rightarrow M_l(Z)$ given by $\lambda (a)(x)=ax$ is
contractive. Similarly, the map
$\rho : D\rightarrow M_r(Z)$ given by $\rho (b)(x)=xb$
is contractive. The maps $\lambda$ and $\rho$ are $w^*$-continuous by \cite[Theorem 4.7.4]{bm}.

Let $u$ be in $M_l(X).$ By Theorem \ref{bim11} for every $\delta
>0$ there exists $y\in Y$ such that $y^{-1}\in Z$,
$\|y\|<1+\delta$ and $\nor{y^{-1}}<1+\delta$. Since the operators
of $M_l(Z)$ and $M_r(Z)$ commute, for all $x\in X$ we have
$$u(y^{-1})yx =\rho (yx)(u(y^{-1})) = u(\rho (yx)y^{-1})=u(y^{-1}yx)=u(x).$$

If $a_\delta =u(y^{-1})y\in C$ then $a_\delta =a $ for all $\delta
>0$ and $\lambda (a)=u.$ It follows that $\lambda $ is surjective. Since
$$\|a\|=\nor{u(y^{-1})y}\leq \|u\|_{M_l(Z)}(1+\delta )^2,
\ \mbox{ for all } \delta > 0,$$ we have that $\|a\|\leq
\|u\|_{M_l(Z)}.$ Thus, $\lambda$ is
isometric.

If $n\in \mathbb{N}$ the  algebras $M_n(C), M_n(D)$ are similar
nest algebras. Repeating the above arguments we can show that
$\lambda $ is $n$-isometric. Hence $\lambda$, and similarly
$\rho,$ is completely isometric. \end{proof}

\begin{example}\label{nestex}
\em{The above result shows that there exists a dual operator space
$Z$ such that $M_l(Z)$ and $M_r(Z)$ not stably isomorphic. Indeed,
from \cite[Example 3.7]{ele2} there exist similar nest algebras
$C$ and $D$ which are not stably isomorphic. The claim now follows
from Theorem \ref{bim12}.}
\end{example}

\end{document}